\documentclass[11pt]{amsart}

\usepackage{amsmath}
\usepackage{amssymb}
\usepackage{amsthm}
\usepackage{graphicx}
\usepackage{subcaption}
\usepackage{hyperref}

\newtheorem{thm}{Theorem}

\newtheorem{prop}[thm]{Proposition}

\newcommand{\E}{\mathbb{E}}

\def \e{\varepsilon}
\newcommand{\be}{\begin{equation}}
\newcommand{\ee}{\end{equation}}

\newcommand{\C}{\mathbb{C}}

\newcommand{\Z}{\mathbb{Z}}

\newcommand{\ol}{\overline}

\newcommand{\ep}{\epsilon}

\renewcommand{\P}{\mathbb{P}}

\theoremstyle{remark} 
\newtheorem{remark}[]{Remark}

\newcommand{\set}[1]{\{#1\}}
\newcommand{\lr}[1]{\left(#1\right)}
\newcommand{\abs}[1]{\left|#1\right|}
\newcommand{\gives}{\rightarrow}
\newcommand{\twiddle}[1]{\widetilde{#1}}
\DeclareMathOperator{\Zeros}{Zeros}
\DeclareMathOperator{\Crit}{Crits}

\title[Random Lemniscate]{The lemniscate tree of a random polynomial}
\author{Michael Epstein, Boris Hanin, and Erik Lundberg}
\date{}
\begin{document}

\begin{abstract}
To each generic complex polynomial $p(z)$
there is associated a labeled binary tree 
(here referred to as a ``lemniscate tree'')
that encodes the topological type of the graph of $|p(z)|$.
The branching structure of the lemniscate tree is determined by 
the configuration (i.e., arrangement in the plane) 
of the singular components of those level sets $|p(z)|=t$
passing through a critical point.

In this paper, we address the question
``How many branches appear in a typical lemniscate tree?''
We answer this question first for a lemniscate tree 
sampled uniformly from the combinatorial class 
and second for the lemniscate tree arising from a random polynomial 
generated by i.i.d. zeros.
From a more general perspective,
these results take a first step toward 
a probabilistic treatment (within a specialized setting)
of Arnold's program of enumerating algebraic Morse functions.
\end{abstract}

\maketitle
%\noindent{\bf No prior NSF support:} Prior work is discussed within the proposal.

\section{Introduction}

Hilbert's sixteenth problem asks for
an investigation of the topology 
of real algebraic curves and hypersurfaces.
An extension of this program,
promoted by V.I. Arnold \cite{Arnold},
is to study the possible equivalence classes
of graphs of generic polynomials
up to diffeomorphism of the domain and range.
Thus, rather than considering a single level set of a polynomial, 
Arnold's problem is concerned with
the whole landscape given by its graph.

The restricted case
to classify graphs arising from the taking the modulus $|p(z)|$
of a generic complex polynomial $p$ 
was solved by F. Catanese and M. Paluszny.\footnote{This problem fits into 
Arnold's setting of 
real polynomials if we equivalently consider 
the square of the modulus and notice that
$p(z) \ol{p(z)}$ has real coefficients as a polynomial in $x$ and $y$.}
They enumerated all possible equivalence classes
by establishing a one-to-one correspondence with the combinatorial class
of labeled, increasing, nonplane, binary trees. 

Motivated by recent studies on the topology of 
\emph{random} real algebraic varieties
\cite{FLL, GaWe0,GaWe2,GaWe3,LLstatistics,Lemni,LundRam,NazarovSodin, NazarovSodin2, SarnakWigman},
it seems natural to investigate 
a probabilistic version of Arnold's problem:
to study the landscape generated by a random polynomial
while focusing on statistics derived from its topological type.
In this paper, 
we investigate the special class mentioned above.
Thus, we consider a random complex polynomial
$p(z)$ and study the induced random binary tree.
We randomize $p$ by sampling independent identically distributed
zeros from a fixed probability measure on the Riemann sphere.

In this particular setting, 
the typical binary trees we observe
(arising from landscapes generated by random polynomials)
do not resemble the ``combinatorial baseline'' provided by sampling 
uniformly from the combinatorial class (see Theorem \ref{thm:poly}
and compare with Theorem \ref{thm:comb}).
Namely, the random tree associated to $p(z)$ typically
has very little branching
(with probability converging to one, 
a shrinking portion of the nodes have two children).

\subsection{Lemniscate trees}
As in \cite{CatPal}, we will call a polynomial $p \in \C[z]$ of degree $n+1$ 
\textit{lemniscate generic} 
(or simply \textit{generic}) if $p'$ has $n$ distinct zeros $w_1, \dots, w_n$ 
such that for each $1 \leq i \leq n$, $p(w_i) \neq 0$ 
and such that $\lvert p(w_i) \rvert = \lvert p(w_j) \rvert$ if 
and only if $i=j$.\footnote{The complement of the set of lemniscate generic polynomials
forms a set of codimension one in the parameter space;
as a result, in many models of random polynomials
(including the ones studied in Section \ref{sec:poly} of this paper)
the condition of being lemniscate generic
holds with probability one.}
To such a polynomial one 
can associate a rooted, nonplane, binary tree $LT(p)$ with $n$ vertices, 
whose vertices are bijectively labeled with the integers from 1 to $n$ such that the labels increase along any path oriented away from the root. We call $LT(p)$ the \textit{lemniscate tree} associated to $p$.  $LT(p)$ encodes the topology of the graph of $|p(z)|$ (or equivalently $|p(z)|^2$ which can be viewed as a Morse function). Its vertices correspond to the $n$ generically distinct critical points of $p$. To construct its edges consider for each critical point $w$ of $p$ the connected component 
\[\Gamma_w\subseteq \set{\abs{p(z)}=\abs{p(w)}}\] 
of the level set $\set{z \in \C: \abs{p(z)}=\abs{p(w)}}$ that contains $w$. The curve $\Gamma_w$, referred to as a \textit{small lemniscate} in \cite{CatPal}, is generically a bouquet of two circles with the self-crossing occuring precisely at $w$. The root vertex of $LT(p)$ corresponds to the critical point with the largest value of $\abs{p}$, which is generically unique. The descendents of the root are defined inductively. The children of a vertex corresponding to a critical point $w$ are the vertices associated to critical points $w'$ for which $\Gamma_{w'}$ is surrounded by one of the petals of $\Gamma_w$ and $\abs{p(w')}$ is the largest possible among all critical points whose singular lemniscates are surrounded by the same petal. Thus, since the critical values of $\abs{p(w)}$ are generically distinct at different critical points, each vertex has zero, one, or two children.
We refer the reader to \cite{CatPal} for more details.

While the lemniscate trees defined above are undirected, it will be convenient for us to adopt a term associated with directed graphs. We may impose an implicit direction on the edges of the trees so that each edge is oriented away from the root (such an  oriented tree is properly called an \textit{out-arborescence} but there will be no confusion here). In this context the number of children a vertex has is its \textit{outdegree}, defined to be the number of directed edges emanating from it. 

\begin{figure}
\centering
  \includegraphics[scale=0.35]{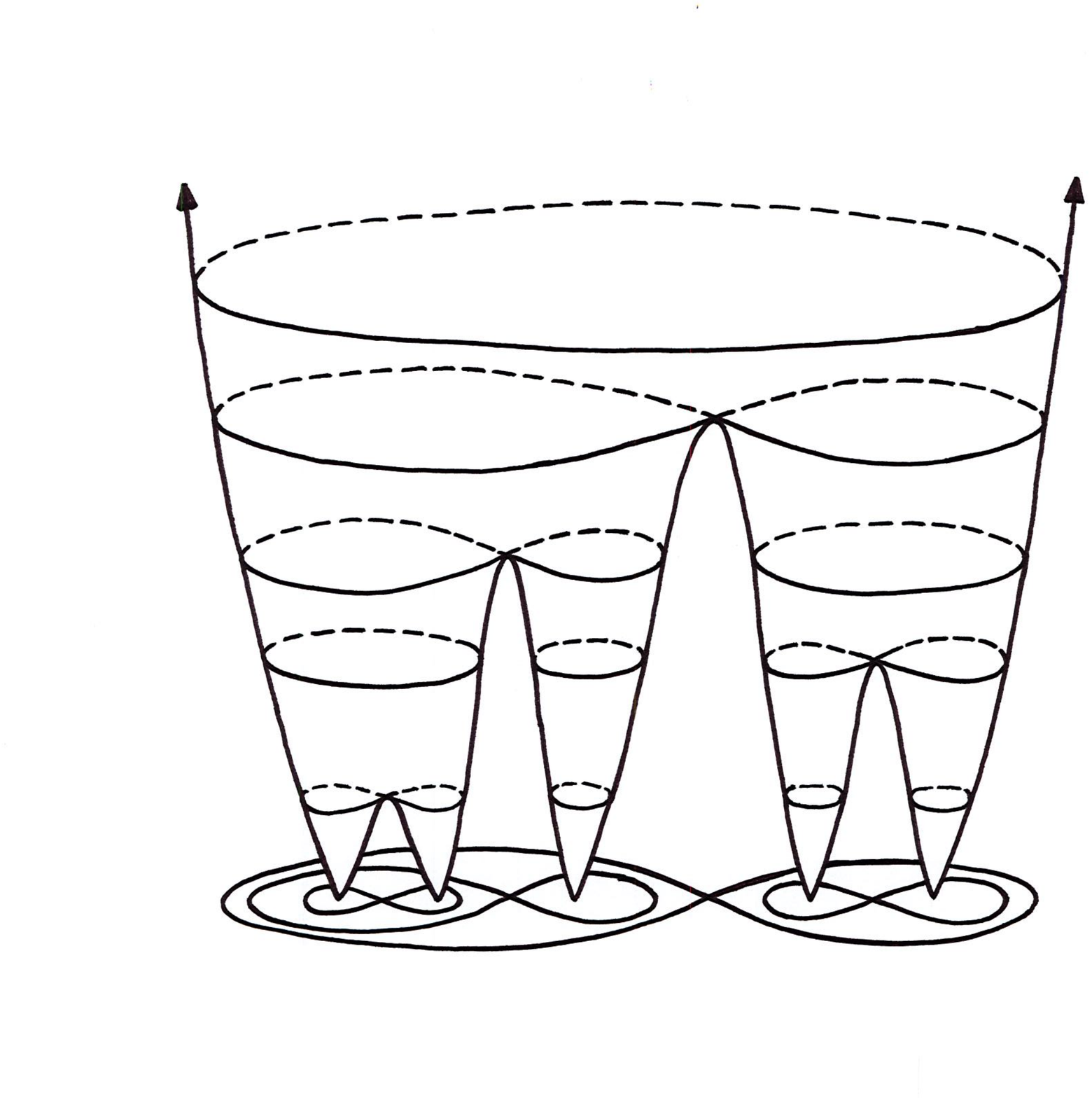}
  \includegraphics[scale=0.45]{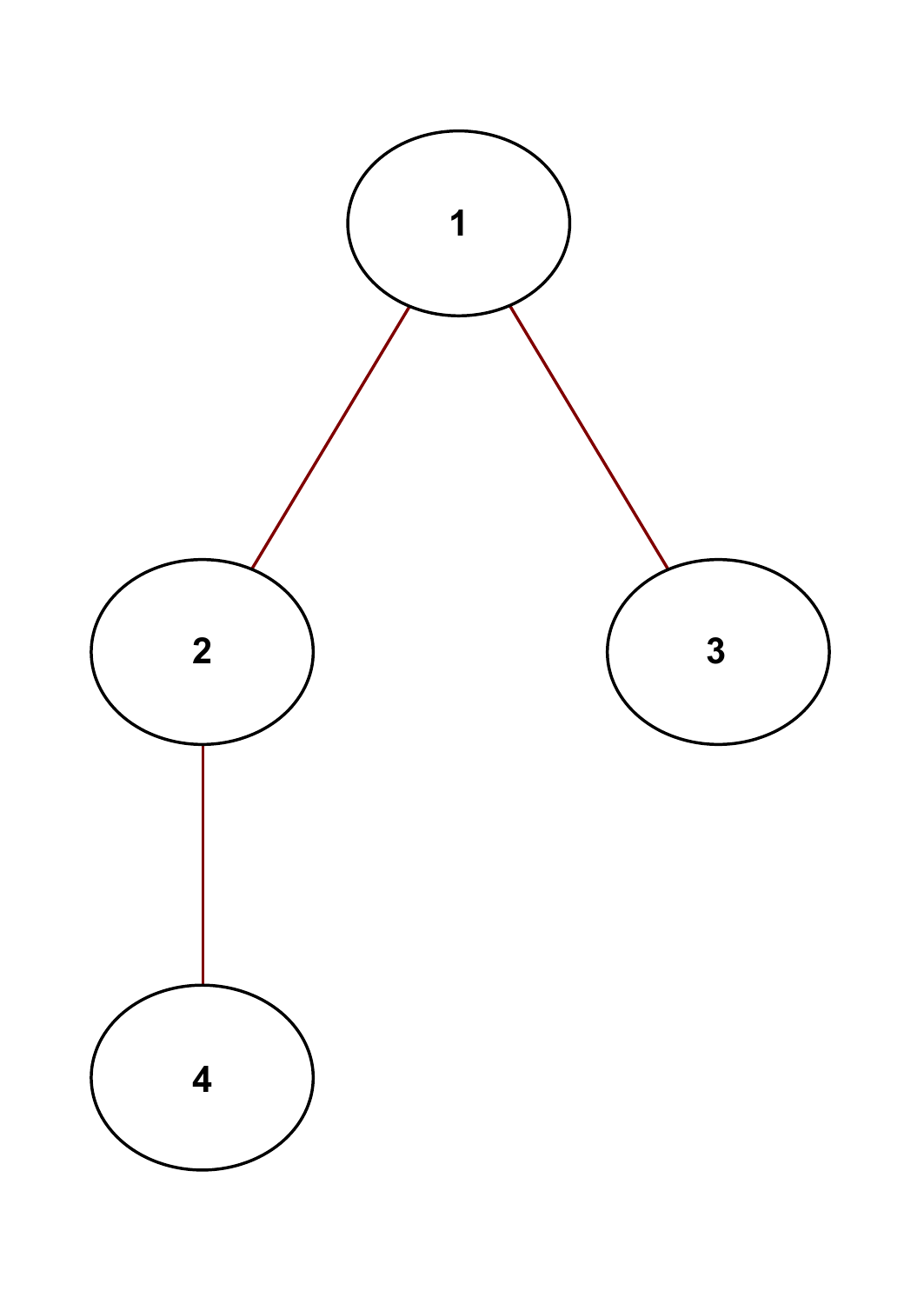}
\caption{Left: an example of a landscape 
(artist's rendition)
generated by a polynomial with five zeros,
along with the projection of the singular component
of each critical level.
Right: the associated lemniscate tree
(each node corresponds to a singular component).
The tree can be constructed using 
the nesting structure of the singular components
along with the ordering of heights of critical values.}
\label{fig:illustration}
\end{figure}

\begin{figure}
\begin{center}
\includegraphics[width=0.3\textwidth]{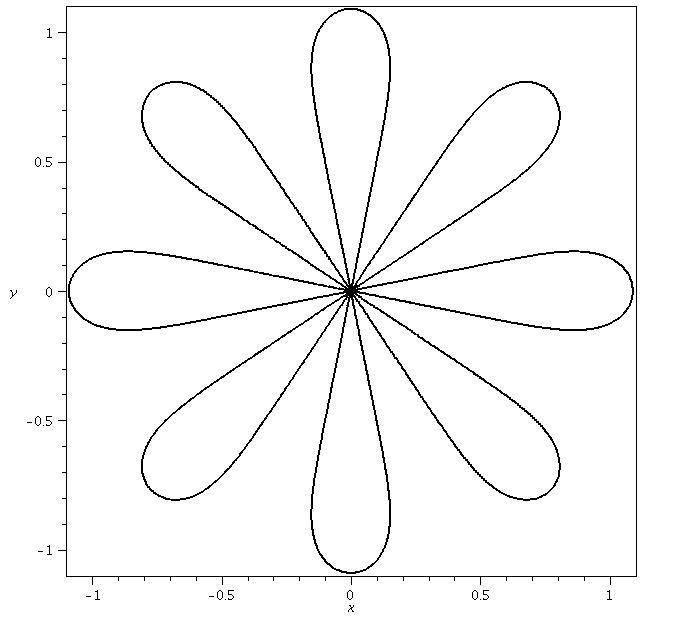}
\caption{The non-generic lemniscate $\{z\in\C:|z^8-1|=1\}$.}
\label{fig:Erdos} 
\end{center}
\end{figure}

\begin{figure}
\centering
\begin{subfigure}{.33\textwidth}
  \centering
  \includegraphics[width=\linewidth]{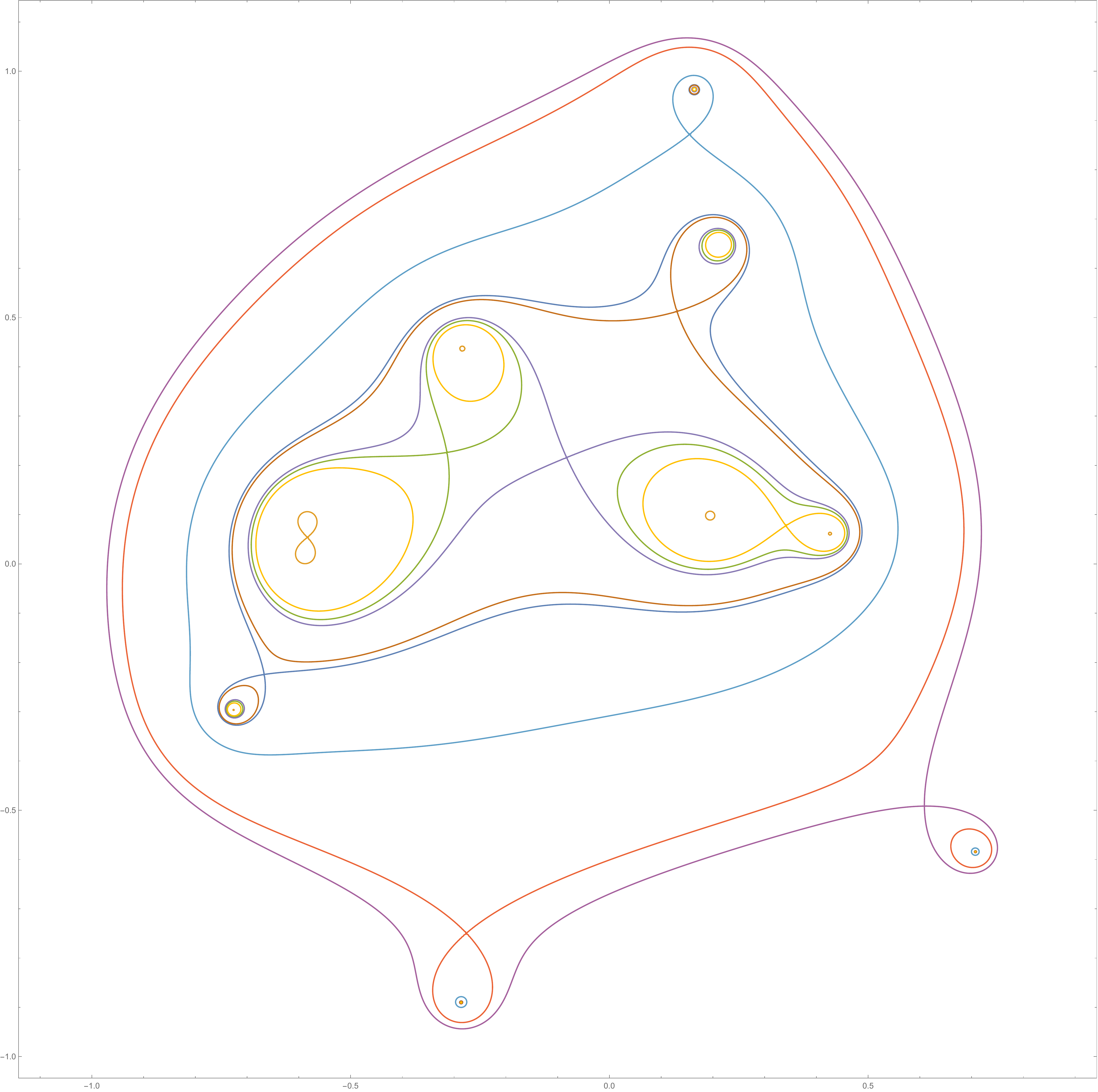}
  \caption{Degree 10}
  \label{fig:sub1}
\end{subfigure}%
\begin{subfigure}{.33\textwidth}
  \centering
  \includegraphics[width=\linewidth]{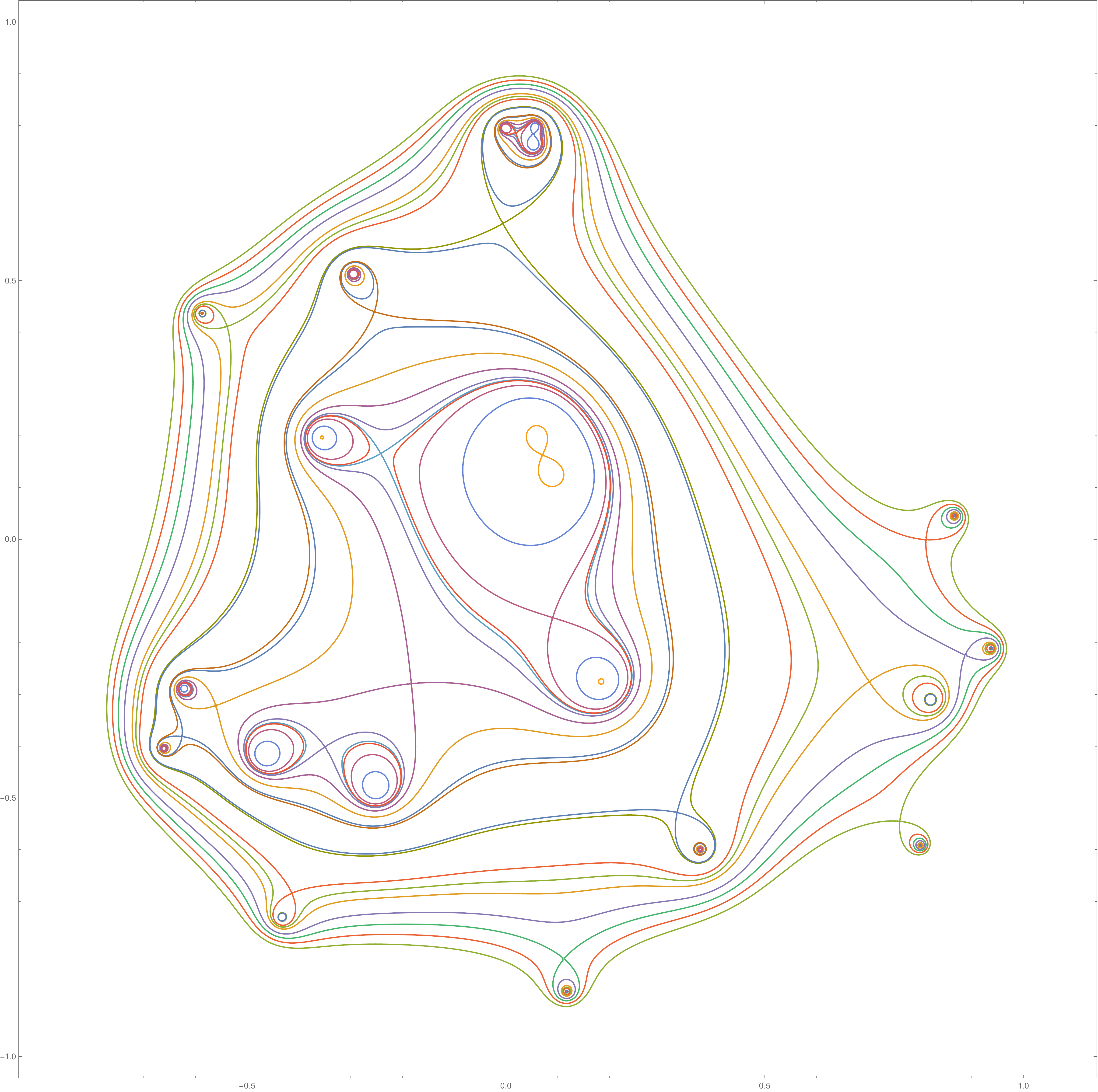}
  \caption{Degree 20}
  \label{fig:sub2}
\end{subfigure}%
\begin{subfigure}{.33\textwidth}
  \centering
  \includegraphics[width=\linewidth]{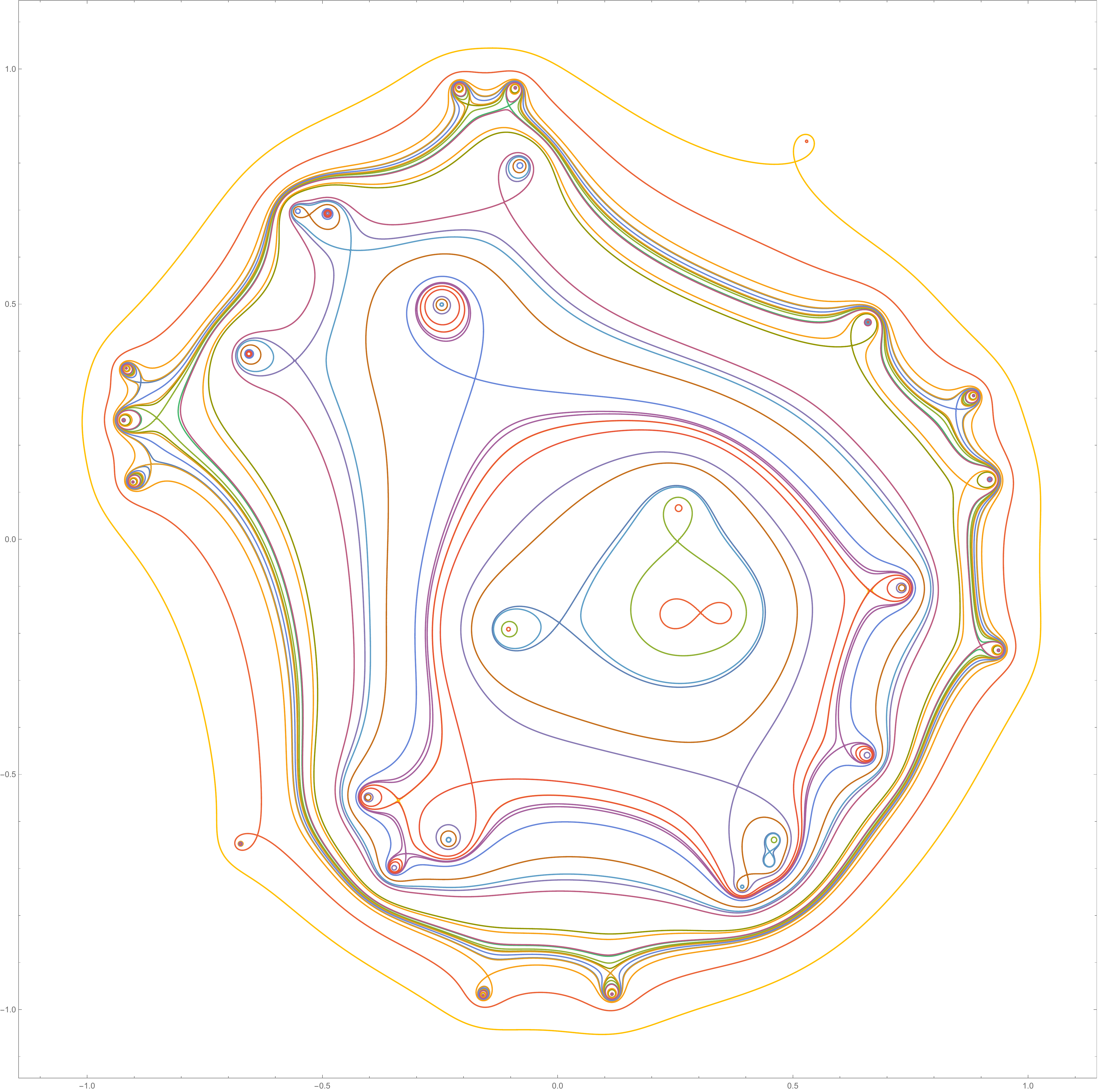}
  \caption{Degree 30}
  \label{fig:sub3}
\end{subfigure}\\
\caption{Lemniscates associated to 
random polynomials generated by sampling i.i.d. zeros
distributed uniformly on the unit disk.
For each of the three polynomials sampled, we have plotted
(using Mathematica)
each of the lemniscates that passes through a critical point.
One observes a trend: most of the singular components have
one large petal (surrounding additional singular components)
and one small petal that
does not surround any singular components. 
Note that only one 
of the connected components in each singular level set is singular
(the rest of the components at that same level are smooth ovals).}
\label{fig:unif}
\end{figure}

\begin{figure}
\begin{center}
\includegraphics[width=0.6\textwidth]{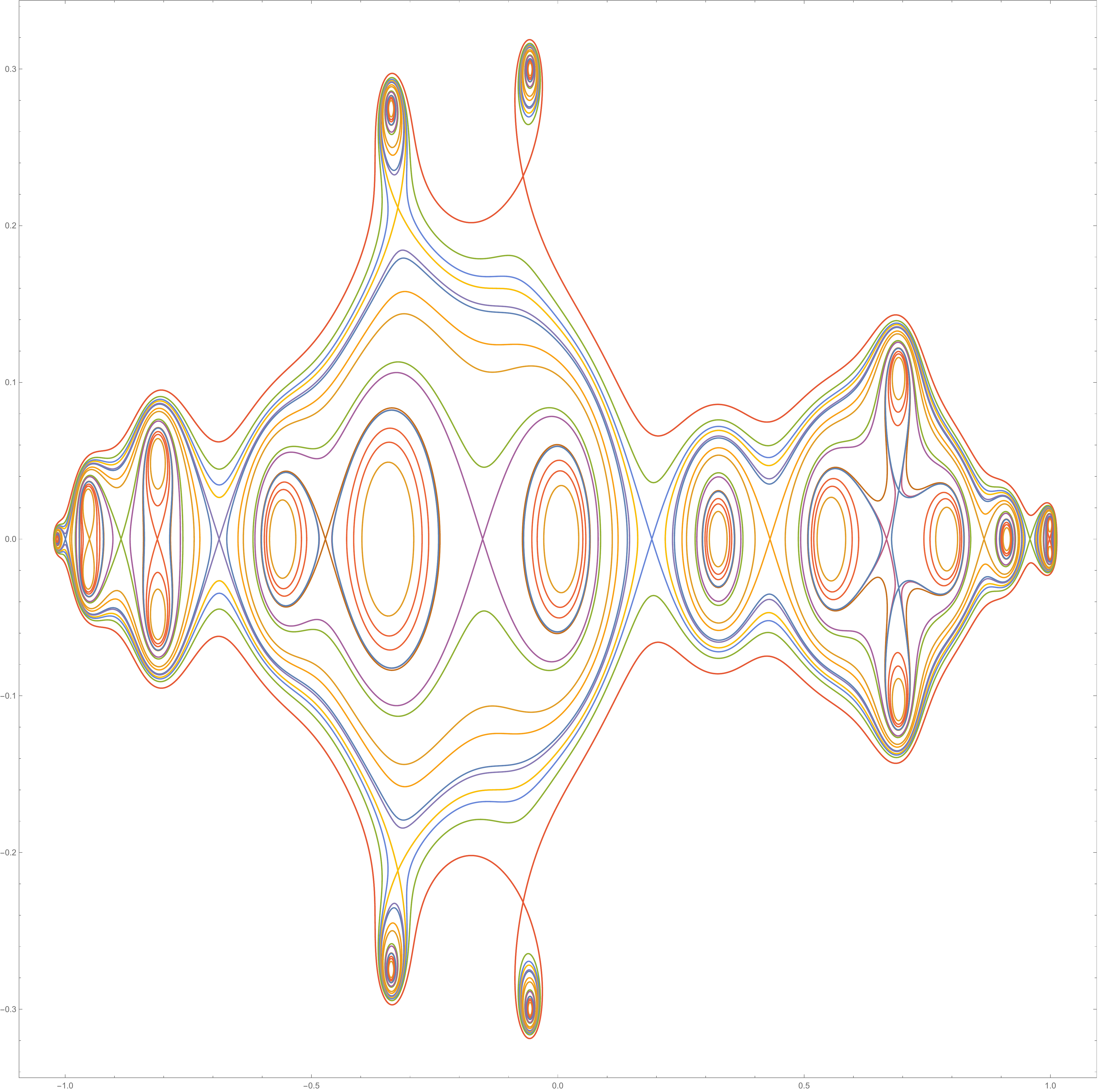}
\caption{Lemniscates associated to a random linear combination
of Chebyshev polynomials with Gaussian coefficients.
Degree $N=20$.  This example is not lemniscate generic
(since we see multiple critical points on a single level set).
However, this model has the interesting feature that it 
seems to generate trees typically having
many branches.  See \S \ref{S:Cheby}.}
\label{fig:Cheb} 
\end{center}
\end{figure}

A simple example illustrating the relation 
between a polynomial and its corresponding tree
is provided in Figure \ref{fig:illustration}. 
The left panel in Figure \ref{fig:illustration} 
displays all the singular level sets 
(each of which may include smooth components in addition to the 
singular component) 
for the modulus of a degree five polynomial,
and the right panel shows the corresponding lemniscate tree. 
To get a sense of what high-degree lemniscates can look like, 
consider Figure \ref{fig:unif},
where the polynomials are generated by sampling i.i.d. zeros
uniform in the unit disk.
Illustrating a highly non-generic lemniscate (for the sake of comparison), 
Figure \ref{fig:Erdos} shows the so-called \emph{Erd\"os lemniscate}
$\{z\in\C:|z^N-1|=1\}$ with $N=8$.
Figure \ref{fig:Cheb} shows 
the singular lemniscates of a random polynomial generated
by a linear combination of Chebyshev polynomials with Gaussian coefficients 
(see \S \ref{S:Cheby} for further discussion of this model).

\subsection{Random lemniscate trees}
In this section we state our main results on the branching in random lemniscate trees. For each $N \geq 1$ we define $LT_N$ to be the set of all lemniscate trees on $N$ vertices. That is, $LT_N$ is the set of all rooted, nonplane, binary trees, with vertices bijectively labeled with the integers from 1 to $N$ such that the labels increase along every path oriented away from the root. As mentioned above, $LT_N$ was shown in \cite{CatPal} to be the space of possible lemniscate trees for generic polynomials of one complex variable. For every $N\geq 1,$ the space $LT_N$ is finite, and our first result concerns the branching structure a tree sampled uniformly at random from $LT_N.$.

\begin{thm}\label{thm:comb}
Let $T_N\in LT_N$ be a lemniscate tree of size $N$ sampled uniformly at random, and let $X_N$ denote the number of vertices of outdegree two in $T_N$. Write $\mu_N$ for its mean and $\sigma_N$ for its standard deviation. Then
$$\mu_N = \left(1-\frac{2}{\pi}\right)N + O(1),$$
and
$$\sigma_N^2 =\left( \frac{4}{\pi^2} + \frac{2}{\pi}-1 \right)N + O(1).$$
Moreover, the rescaled random variable $\sigma_N^{-1} \left( X_N - \mu_N \right)$ converges in distribution to a standard Gaussian random variable as $N\gives\infty.$ 
\end{thm}

We note that the asymptotic for the mean follows from
the asymptotic for the mean number of leaves which was computed recently in \cite{Bona1} (cf. \cite{Bona2}, \cite{BonaPittel}, \cite{Bona2014}),
where the same class of trees was referred to as $1$-$2$ trees.

By a standard application of Chebyshev's inequality, we see that Theorem \ref{thm:comb} implies that the number of nodes of outdegree two is concentrated about its mean. Indeed, choosing $0<\alpha<1/2$ we have
$$P(|X_n - \mu_n| > n^{\alpha+1/2} ) \leq \frac{\sigma_n^2}{n^{\alpha+1/2}} = O(n^{\alpha-1/2}) = o(1), \quad \text{as } n \rightarrow \infty.$$
Therefore, for a uniformly randomly sampled $T_N\in LT_N$, one expects a constant proporition of its vertices to have two children. Our next result concerns the number of outdegree $2$ nodes in the lemniscate tree of a random polynomial. Formally, we equip the space of polynomials of degree $N$ with a measure under which zeros are chosen i.i.d. on the Riemann sphere, and push forward this measure to $LT_{N-1}$ under the map that associates to a generic polynomial its lemniscate tree. 

 \begin{thm}\label{thm:poly}
 Let $p_N$ be a random polynomial of degree $N$ whose zeros are drawn 
 i.i.d. from a fixed probability measure $\mu$ on $S^2$ that has a 
 bounded density with respect to the uniform (Haar) measure. 
 Then for every $\ep>0$ there exists $C_\ep$ so that the number 
 $Y_N$ of nodes of outdegree two in the lemniscate tree 
 associated to $p_N$ satisfies
 \[\E Y_N \leq C_\ep N^{\frac{1}{2}+\ep}.\]
 \end{thm}
Although Theorem \ref{thm:poly} does not give variance estimates and asymptotic normality as in Theorem \ref{thm:comb}, it does show that in contrast to sampling uniformly from $LT_N$ the lemniscate tree of a random polynomial has almost all nodes with outdegree at most one. It also provides a weak concentration inequality for the random variables $Y_N$. 
Namely, for any $\varepsilon_0 > 0$
$$P \left( Y_N > N^{\frac{1}{2}+ \e_0} \right) \leq \frac{\E Y_N}{N^{\frac{1}{2} + \e_0}} \rightarrow 0, \quad \text{as } N \rightarrow \infty.$$

This sparse branching for the lemnisate trees of a random polynomials
is closely related to the \emph{pairing of zeros and critical points}
for random polynomials studied by the second author \cite{Hanin1, Hanin2, Hanin} and taken up in \cite{OW} as well. These articles roughly show that for the random polynomials we consider, each zero of $p$ has, with high probability, a paired critical point in its $1/N$ neighborhood. As we show in the proof of Theorem \ref{thm:poly}, when such a pairing occurs, the singular lemniscate $\Gamma_w=\set{\abs{p(z)}=\abs{p(w)}}$ passing through the critical point $w$ of $p$ that is paired to a zero $z$ is likely to have a small petal surrounding $z$ and no other zeros (and hence no other singular lemniscates either), causing the corresponding vertex to have outdegree (at most) $1.$ 

A useful heuristic for understanding the critical point pairing (which we combined with 
a topological argument in order to prove Theorem \ref{thm:poly}) is in terms of electrostatics on $S^2$, where critical points of $p$ are viewed as equilibria of the field generated by a logarithmic potential with positively charged point particles at the zeros and negatively charged particles at the poles, counted with multiplicity. 

The contribution to the electric field from the high-order pole at infinity (which is best understood after changing coordinates by $z = 1/\zeta $) 
is balanced by the electric field from an individual zero
in a neighborhood with radius of order $1/N$.
In this neighborhood, the additional influence of other zeros of $p$ is typically of lower order, causing an almost deterministic pairing of zeros and critical points. 
We refer the reader to \S \ref{sec:poly} below and to \cite[\S 1]{Hanin} for more details. 
 
%This leads to the following interesting problem.

%\vspace{0.1in}

%\noindent {\bf Open problem:} Investigate the sharp concentration around the mean for the random variable $\Lambda_n.$ In other words, what is the best estimate one can get for $\mathbb{P}\{|\Lambda_n - \EE |\Lambda_n|| > t\}$?

\begin{remark}
	There is a fair amount of ``universality'' expressed in Theorem \ref{thm:poly} in that the distribution $\mu$ is rather arbitrary.
	What if the polynomial is instead sampled
	using random coefficients in front of some choice of basis?
	Based on simulations, 
	the lemniscate trees again seem to have a shrinking portion of
	nodes with two children in a wide variety of such models including
	most of the well-studied Gaussian models
	(the Kostlan model, the Weyl model, the Kac model).
	In fact, the only exception we observed
	was a model based on Chebyshev polynomials 
	(see the empirical evidence presented below in the last section).
	In another (more exotic) direction, one may consider
	randomizing the construction of polynomial ``fireworks''
	described in \cite[\S 4]{FKV}
	in order to produce polynomials whose trees have many branches.
\end{remark}

\begin{remark}
	As a future direction of study it seems natural to investigate random \emph{rational} functions on the Riemann sphere.	
	A combinatorial scheme for classifying associated topological types
	was developed in \cite{Catanese2}.
	What positive statements can one make on the typical topological type?
	The results in \cite{Lemni},
	investigating a fixed level set of a random rational function 
	(defined as the ratio of two random polynomials from the Kostlan ensemble), may lead to some insight in this direction.
	However, we generally anticipate the case of rational functions to have a much different flavor than the case of polynomials; 
	not only is the underlying combinatorial class more complicated, but there is no longer a ``polarization'' caused by having a high-order pole at infinity.
\end{remark}

\begin{remark}
	Another natural direction of study,
	returning to Arnold's problem mentioned at the beginning of the introduction,
	would be to investigate the topological type of a random homogeneous polynomial in projective space.
	The underlying classification problem in this case is still unsolved;
	L. Nicolaescu classified generic Morse functions on the
	$2$-sphere \cite{Nicolaescu} 
	and enumerated them in terms of their number of critical points,
	but it is not known which types can be realized within each
	space of polynomials of given degree (and even less is known in more than two variables) \cite{Arnold}.
	At this stage, we suggest investigating a coarser structure,
	such as the so-called ``merge tree'' \cite[\S VII.1]{EdelHarer},	
	associated to the graph of a random real homogeneous polynomial
	of degree $d$ in $n+1$ variables
	(while pursuing asymptotic estimates as $d \rightarrow \infty$ for statistics defined on the merge tree).
\end{remark}

\subsection{Outline of the paper}

The Gaussian limit law stated in
Theorem \ref{thm:comb}
will be established using \emph{perturbed singularity analysis},
a method from analytic combinatorics.
Specifically, in \S \ref{S:comb-pf}, we will apply a
result from \cite{Flajolet}
to a bivariate generating function 
that was derived in \cite{Appendix}.
We prove Theorem \ref{thm:poly}
in Section \ref{sec:poly}
by establishing a prevalence of
small lemniscate petals
adapting the method from \cite{Hanin}
for studying pairing between zeros and critical points
of random polynomials.
In Section \ref{S:Cheby},
we present some empirical results concerning
a certain model of random polynomials for which 
the lemniscate trees appear to have 
on average asymptotically one third of their nodes
being of outdegree two.

\section{Sampling uniformly from the combinatorial class: proof of Theorem \ref{thm:comb}}\label{S:comb-pf}

Let $a_{n,k}$ denote the number of lemniscate trees 
of size $n+1$ with $k$ nodes of outdegree two,
and consider the bivariate generating function
$$F(z,u) = \sum_{n,k \geq 0 } \frac{a_{n,k}}{n!} u^k z^n.$$

In \cite{Appendix}, an explicit formula for the function
$F(z,u)$ is derived by showing that $F$ satisfies a 
first-order PDE that can be solved 
explicitly using the method of characteristics.
This results in the following analytic description in 
terms of elementary functions
\begin{equation}\label{eq:bivariate}
F(z,u) =  \left[ \cosh\left(\frac{z}{2}\sqrt{1-2u}\right) - \frac{\sinh\left(\frac{z}{2}\sqrt{1-2u}\right)}{\sqrt{1-2u}} \right]^{-2}.
\end{equation}

There is a well-established theory for deriving probabilistic results from
bivariate generating functions such as $F(z,u)$.
For a detailed overview, see the authoritative text \cite[Ch. IX]{Flajolet} 
by Ph. Flajolet and R. Sedgewick;
here we briefly review the connection in the current context.
The basic link is that we arrive at the so-called
\emph{probability generating function}
by considering a normalized coefficient extraction involving $F(z,u)$.
Namely, using $[z^n]$ to denote the operation of extracting 
the $z^n$-coefficient, the univariate polynomial in $u$, 
given by
$$p_n (u) = \frac{[z^n ]F(z, u)}{[z^n ]F(z, 1)},$$
is the probability generating function for the random variable $X_n$ 
defined (as in the statement of Theorem \ref{thm:comb}) 
as the number of nodes of outdegree two
in a random lemniscate tree of size $n$.
That is, if a lemniscate tree of size $n$ is sampled uniformly
at random, the probability that it has $k$ nodes of outdegree two
is given by the coefficient of $u^k$ in $p_n(u)$.
From this, one can easily compute the mean and 
variance using simple operations.
Furthermore, a more detailed
complex analysis of the singularity structure of 
bivariate generating functions such as $F(z,u)$
can be used to establish probabilistic limit laws.  

Concerning the case at hand, viewing $u$ as a complex parameter,
the function $F(z,u)$ is amenable to perturbed singularity analysis
and falls under the ``movable singularities schema''
described in \cite{Flajolet};
as $u$ varies in a neighborhood of $u=1$,
the location of the (nearest to the origin) singularity of $F(z,u)$ 
moves while the nature of this singularity is preserved.
This allows us to establish a Gaussian limit law by apply the following result
restated from \cite[Thm. IX.12]{Flajolet}.
\begin{thm}\label{thm:Flajolet}
Let $F(z,u)$ 
be a function that is
bivariate analytic at 
$(z, u) = (0,1)$ 
and has non-negative coefficients. 
Assume the following conditions hold:
\begin{itemize}
\item[(i)] Analytic perturbation: there exist three functions $A, B, C,$ 
analytic in a domain 
$\mathcal{D} = \{ |z| \leq r \} \times \{ |u - 1| < \e \}$, 
such that the following representation holds
in some neighborhood of $(0,0)$, 
with $\alpha \notin \Z_{\leq 0}$,
$$F(z, u) = A(z, u) + B(z, u)C(z, u)^{-\alpha}.$$
Furthermore, in $|z| \leq r$, 
there exists a unique root $\rho_1$ of the
equation $C(z, 1) = 0$, this root is simple, 
and $B(\rho_1, 1) \neq 0$.
\item[(ii)] Non-degeneracy: one has 
$\partial_z C(\rho_1, 1) \cdot \partial_u C(\rho_1, 1) \neq 0$, 
ensuring the existence of a non-constant analytic 
function $\rho(u)$ near $u = 1$, 
such that $C(\rho(u), u) = 0$ and $\rho(1) = \rho_1$.
\item[(iii)] Variability: one has
$$v(\beta):=\frac{\beta''(1)}{\beta(1)} + \frac{\beta'(1)}{\beta(1)} - \left(\frac{\beta'(1)}{\beta(1)}\right)^2 \neq 0,$$
where $\beta(u) = \rho(1) \rho(u)^{-1}$.
\end{itemize}
Then, the random variable with 
probability generating function
$$p_n (u) = \frac{[z^n ]F(z, u)}{[z^n ]F(z, 1)}$$
converges in distribution (after standardization) 
to a Gaussian random variable with a speed of convergence $O(n^{-1/2})$. 
\end{thm}

\subsection*{Verifying condition (i)}
Let $G(z,u) = \cosh\left(\frac{z}{2}\sqrt{1-2u} \right) - \frac{\sinh\left(\frac{z}{2}\sqrt{1-2u}\right)}{\sqrt{1-2u}}$ 
so that we have $F(z,u) = G(z,u)^{-2}$.
Note that $G(z,u)$ is an entire function of $z$ for each fixed $u \in \C$,
and is non-constant for $u \neq \frac{1}{2}$.
Thus, for $u \neq \frac{1}{2}$, 
$F(z,u)$ is meromorphic with 
poles at the zeros of $G(z,u)$ and no other singularities.
First setting $u=1$, 
we find that the zeros of $G(z,1)$
are at $z=\frac{\pi}{2} + 2 \pi i k, k \in \Z$.
Among these, $\rho_1 = \frac{\pi}{2}$ is nearest to the origin.
We compute $\partial_z G(\rho_1,1) = \frac{-\sqrt{2}}{2} \neq 0$,
which shows that $\rho_1$ is a simple root.
This completes the verification of condition (i) in Theorem \ref{thm:Flajolet},
where $A = 1, B = 1$ are taken to be constant, $\alpha = 2$,
and $C(z,u) = G(z,u)$.

\subsection*{Verifying condition (ii)}
Having shown above that $\partial_z G(\rho_1,1) \neq 0$,
we only need to check that $\partial_u G(\rho_1, 1) \neq 0$.
We find $\partial_u G(\rho_1, 1) = \frac{\sqrt{2}}{4}\left(2-\pi \right)$.
This verifies condition (ii),
where the desired function $\rho(u)$
is guaranteed to exist by the implicit function theorem.
Furthermore, 
we can describe $\rho(u)$ explicitly by solving for $z$ in $G(z,u)=0$:
\[
		\rho(u) = \frac{2}{\sqrt{1-2u}}\tanh^{-1}\left(\sqrt{1-2u}\right) = \frac{1}{\sqrt{1-2u}}\log\left(\frac{1+\sqrt{1-2u}}{1-\sqrt{1-2u}}\right),
\]
where we choose the principal branch for the logarithm
so that $\rho(1)=\rho_1 = \pi/2$.
For $|u-1|$
sufficiently small, $\sqrt{1-2u}$ is near $i$,
ensuring analyticity of the function $\rho(u)$.

\subsection*{Verifying condition (iii)}

Let $\beta(u) = \rho(1)/\rho(u)$. 
Then 
$$\beta'(u) =-\frac{\pi}{2}\rho(u)^{-2}\rho'(u),$$ and 
$$\beta''(u) 
=-\frac{\pi}{2}\left( -2\rho(u)^{-3}(\rho'(x))^2 +\rho(x)^{-2}\rho''(x) \right).$$
Thus,
\begin{align*}
	v(\beta) &= \frac{\beta''(1)}{\beta(1)} + \frac{\beta'(1)}{\beta(1)}- \left(\frac{\beta'(1)}{\beta(1)}\right)^2\\
		    &= \left(\frac{8}{\pi^2}-1\right) + \left(1-\frac{2}{\pi}\right) - \left(1-\frac{2}{\pi}\right)^2\\
		    &=\frac{4}{\pi^2} + \frac{2}{\pi} -1 \neq 0.
\end{align*}

\noindent
We conclude that Theorem \ref{thm:Flajolet} applies, 
and the random variable
$N_2(T_n)$ after rescaling converges in distribution 
to a Gaussian variable with a speed of convergence O$(n^{-1/2})$. 

\subsection*{Mean and variance}
As pointed out in \cite{Flajolet} in the remarks after the proof of 
Theorem IX.12, 
the mean $\mu_n$ and variance $\sigma_n^2$ are given by 
$$\mu_n = \frac{\beta'(1)}{\beta(1)} n +\text{O}(1) = \left(1-\frac{2}{\pi}\right)n + \text{O}(1),$$
and 
$$\sigma_n^2 = v(\beta) \cdot n + \text{O}(1) = \left(\frac{4}{\pi^2} +\frac{2}{\pi} -1 \right)n +\text{O}(1).$$

\section{Proof of Theorem \ref{thm:poly}}\label{sec:poly}
\noindent For a polynomial $p$ of one variable we define 
\[\Zeros(p):=\set{\xi \in S^2\text{ s.t. } p(\xi)=0},\qquad\Crit(p):=\set{w\in S^2\text{ s.t. } dp(w)=0}.\] 
Instead of working the ``usual'' holomorphic coordinates $S^2\backslash\set{\infty}\gives \C,$ it will be more convenient to perform our computations in coordinates $S^2\backslash\set{0}\gives \C$ centered at the point at infinity. 
That is, we write 
\[p_N(w):=\frac{1}{w^N}\prod_{j=1}^N \lr{w-\xi_j},\]
where $\xi_j$ are drawn i.i.d. from $\mu.$ 
Let us emphasize that whenever a condition like $\abs{\xi}\leq N^{\Delta}$ appears below for some $\xi\in S^2$, the quantity $\abs{\xi}$ is computed in this system of coordinates. In particular, denoting by $\zeta=1/\xi$ the image of $\xi$ in the usual coordinates centered at $0$, 
our condition $\abs{\xi}\leq N^{\Delta}$ is the same as $\abs{\zeta}\geq N^{-\Delta}.$ Associated to each $w\in \Crit(p_N)$ 
is the singular component $\Gamma_w$ of the lemniscate
\[ \Lambda_w := \set{z\in S^2 \text{ s.t. } \abs{p_N(z)}=\abs{p_N(w)}},\]
that passes through $w$. That is, among the connected components of $\Lambda_w$, we define $\Gamma_w$ to be the one that contains $w$. For a \emph{generic} polynomial (a condition that holds with probability one in our model), there are $N-1$ distinct singular lemniscates (one passing through each critical point), each having a unique singular component that is topologically a bouquet of two circles. We call these two circles \textit{the petals} of $\Gamma_w.$ For the arguments below we fix an auxiliary parameter $r\gg 1$ such that
\begin{equation}
\arg\lr{1+\frac{1}{r}e^{i\theta}}\in \lr{-\frac{1}{10},\frac{1}{10}},\qquad \forall \theta\in [0,2\pi].\label{E:r-def}
\end{equation}
We study the behavior of the lemniscate tree of $p_N$ by considering for each $\xi\in \Zeros(p_N)$ 
the event
\[S_{\xi,N}:=\left\{\exists ! \, w \in \Crit(p_N)~\bigg|~
\substack{\abs{\xi-w}<\frac{r\abs{\xi}}{N} ~~\text{and at least one petal of }\Gamma_{w}\\ \text{is contained in the disk of radius }\frac{4r\abs{\xi}}{N}\text{centered at }\xi}\right\}.\]
When the event $S_{\xi,N}$ occurs,
we will say that \textit{$\Gamma_{w}$ has a small petal surrounding $\xi$},
and we refer to $w$ as the \textit{paired critical point} of $\xi$. 
We also consider the events
\[B_{\xi,N}:=\left\{\abs{\xi-\xi'}>\frac{4r\abs{\xi}}{ N},\qquad \forall \xi'\in \Zeros(p_N)\backslash\set{\xi}\right\}.\]
To prove Theorem \ref{thm:poly}, we begin by observing that
\begin{align}
\label{E:lb}\#\set{\text{outdegree at most }1\text{ nodes in lemniscate tree of }p_N}\geq \sum_{\xi\in \Zeros(p_N)}    {\bf 1}_{S_{\xi,N}\cap B_{\xi,N}},
\end{align}
where ${\bf 1}_S$ denotes the indicator function of the event $S.$ To see that \eqref{E:lb} holds, observe that if the events $S_{N,\xi}\cap B_{N,\xi}$ and $S_{N,\xi'}\cap B_{N,\xi'}$ occur for some zeros $\xi\neq \xi'$, then the corresponding paired critical points are also distinct since the spacing of zeros ensured by $B_{\xi,N}$ and $B_{\xi',N}$ is larger than the sum of the distances between the zeros to their paired critical points given by $S_{\xi,N}$ and $S_{\xi',N}$. 
Moreover, when the event $S_{\xi, N}\cap B_{\xi,N}$ occurs, the vertex in the lemniscate tree of $p_N$ that corresponds to the critical point $w$ paired to $\xi$ has outdegree at most $1$.
Indeed, one of its petals surrounds only one zero, namely $\xi$,
and the argument principle then implies that $p_N$ maps the interior of that petal
univalently to a disk (with radius given by $|p_N(w)|$).
This implies that there are no critical points of $p$ in the interior of the petal
(i.e., the vertex in the lemniscate tree of $p_N$ that corresponds to $w$ has outdegree at most $1$).

%Next, consider a critical point $w\in \Crit(p_N)$ for which the corresponding vertex in the lemniscate tree $LT(p_N)$ of $p_N$ has outdegree $2.$ Each of the petals of $\Gamma_w$ must surround the singular component of the lemniscate passing through some element of $ \Crit(p_N)$. We denote these critical points by $w_1,w_2.$ For $j=1,2$, the curves $\Gamma_{w_j}$ must surround at least two zeros of $p_N.$ Indeed, the number of zeros surrounded by any connected component of the level set of the modulus of a holomorphic function is one larger than the number of surrounded critical points (Theorem \cite[Thm ?]{?}). 
This proves \eqref{E:lb} and shows that
\begin{align}
\label{E:lb-avg}\E{\#\set{\text{vertices in }LT(p_N)\text{ with at most one child}}}~~\geq~~ N \cdot \P \lr {S_{\xi,N}\cap B_{\xi,N}}.
\end{align}
To obtain a lower bound for the probability of $S_{\xi,N}\cap B_{\xi,N}$, note that for any $\Delta\in [0,1/4)$ and any $\xi\in \Zeros(p_N)$ we have
\[\mathbb P\lr{\abs{\xi}\leq N^\Delta} = 1 + O\lr{N^{-2\Delta}}\]
since the measure $\mu$ assigns to a ball of radius $N^{-\Delta}$ centered at any point on $S^2$ (in particular at $0$) a mass on the order of its volume. Using that $\mu$ has a bounded density with respect to the uniform measure on $S^2,$ we have
\begin{align*}
  \mathbb P\lr{B_{\xi,N}}&= \mathbb P\lr{B_{\xi,N}~|~\abs{\xi}\leq N^\Delta}\lr{1+O(N^{-2\Delta})}\\
&=\left[1-(N-1)\int_{\abs{\xi}\leq N^\Delta} \mu\lr{\left\{\zeta ~\bigg|~\abs{\xi-\zeta}< \frac{r\abs{\xi}}{N}}\right\} d\mu(\xi)\right]\lr{1+O(N^{-2\Delta})}\\
&=\left[1-O\lr{N^{-1+2\Delta}}\right]\lr{1+O(N^{-2\Delta})}\\
&= 1 + O\lr{N^{-2\Delta}} + O(N^{-1+2\Delta}).
\end{align*}
 Therefore, since for $\Delta\in [0,1/4),$ we have $-2\Delta>-1 + 2\Delta,$ we find that
\begin{align}
\P\lr{S_{\xi,N}\cap B_{\xi,N}}&=   \P\lr{S_{\xi,N}}+ O(N^{-2\Delta})\geq \int_{\abs{\xi}\leq N^{\Delta}} \mathbb P\lr{S_{\xi,N} ~|~  \xi} d\mu(\xi) + O(N^{-2\Delta})\label{E:sn-est},
\end{align}
where the notation in the last integral is that we've conditioned on the position of $\xi.$ We now fix $\Delta\in [0,1/4)$, a deterministic sequence $\xi=\xi_N$ with $\abs{\xi}\leq N^{\Delta}$,
and consider the random polynomials
\[p_{\xi,N}(w):=\frac{1}{w^N}\lr{w-\xi}\prod_{j=1}^{N-1}\lr{w-\xi_j},\]
conditioned to have a zero at $\xi$ and with $\xi_j$ drawn i.i.d. from $\mu$ for $j=1,\ldots, N-1.$ We slightly abuse notation and continue to write $S_{\xi,N}$ for the event that (the fixed zero) $\xi$ has a paired critical point $w_{\xi,N}$ with a small petal surrounding $\xi$, so the conditional probability appearing in the integrand
in (\ref{E:sn-est}) is henceforth simply denoted as $\mathbb P\lr{S_{\xi,N}}$. Theorem \ref{thm:poly} follows from \eqref{E:sn-est} once we show that there exists $N_0\geq 1$ and $C_{\Delta}>0$ so that for all $N\geq N_0$
\begin{equation}\label{E:Goal1}
\inf_{\abs{\xi}\leq N^{\Delta}}\mathbb P\lr{S_{\xi,N}}\geq 1- C_{\Delta} N^{-2\Delta}.
\end{equation}
To show \eqref{E:Goal1}, we revisit the proof of the main theorem in \cite{Hanin}. To state the precise estimate we will use, we set some notation. Critical points of $p_{N,\xi}$ are solutions to $E_N(w) = 0$, where
\begin{equation}\label{eq:electric}
E_N(w) = d\log p_{N,\xi}(w)= -\frac{N}{w} + \frac{1}{w-\xi}+\sum_{j=1}^{N-1} \frac{1}{w-\xi_j}.
\end{equation}
As in \cite[\S 4]{Hanin}, observe that
\begin{equation}\label{eq:expectedfield}
\E{E_N(w)} = -\frac{N}{w} + \frac{1}{w-\xi} + (N-1)\int_{\C}\frac{d\mu(z)}{w-z},
\end{equation}
and
\begin{equation}
\twiddle{E}_N(w) := E_N(w) - \E{E_N(w)} = \sum_{j=1}^{N-1} \frac{1}{w-\xi_j}.
\end{equation}
In the computations below, the Cauchy-Stieltjes transform $\int_{\C}\frac{d\mu(z)}{w-z}$ appearing in \eqref{eq:expectedfield} plays no significant role (it only shift the locations of critcal points in a deterministic way so that \eqref{E:w-def} below has an additional deterministic $1/N$ correction). Hence, we will assume that it is identically $0$ (i.e. we reduce to the case when $\mu$ is the uniform measure on $S^2$). For each $\Delta \in [0,1/4)$ and all $\xi$ with $\abs{\xi}\leq N^{\Delta}$, the average critical point equation $\E{E_N(w)}=0$ has a unique solution
\begin{equation}
w_{\xi,N} := \xi \left( 1 - \frac{1}{N} \right)^{-1}\label{E:w-def}
\end{equation}
near $\xi$. Note that
\[\abs{w_{\xi,N}-\xi}= \frac{\abs{\xi}}{ N-1}.\]
We will argue in \S \ref{S:main-prop-pf} below that the technique in \cite{Hanin} gives the following proposition. 
\begin{prop}\label{prop:main}
Fix $\Delta \in (0,1/4).$ For each $\xi$ let $D_{N,\xi}$ denote the disk of radius $4r \abs{\xi}/N $ centered at $\xi$. There exists $\gamma>0$ and a constant $C_{\Delta} $ so that the event
\[X_{N,\xi,\Delta}=\left\{\sup_{w \in D_{N,\xi}} \left| \twiddle{E}_N(w) \right| \leq N^{1-\gamma}/\abs{\xi}\right\}\]
occurs with high probability:
\begin{equation}\label{eq:disk}
\inf_{\abs{\xi}\leq N^{\Delta}}\P \left( X_{N,\xi,\Delta} \right) \geq 1-C_{\Delta}\cdot N^{-2\Delta}.
\end{equation} 
\end{prop}
Assuming Proposition \ref{prop:main} for the moment, we complete the proof of \eqref{E:Goal1} and hence of Theorem \ref{thm:poly} by showing that the event $X_{N,\xi,\Delta}$ (or more precisely $X_{N,\xi,\Delta}\cap B_{N,\xi}$), whose probability is estimated in \eqref{eq:disk}, is contained in the event $S_{\xi,N}$. Suppose that $X_{N,\xi,\Delta}$ occurs. Then, as in \cite{Hanin}, by Rouch\'e's Theorem applied to $E_N$, there exists $N_0$ so that for all $N\geq N_0$ 
there is a unique $w \in \Crit(p_N)$ satisfying
\[ \abs{w - \xi} < r\abs{\xi}/N\]
with probability at least $1-C_\Delta\cdot N^{-2\Delta}.$ 
Indeed, write $\Gamma$ for the boundary of the disk of radius $r\abs{\xi}/N$ centered at $\xi.$ Since $r\gg 1,$ 
the curve $\Gamma$ winds around $w_{N,\xi}$ (defined in \eqref{E:w-def}) for all $N.$ Moreover, by the triangle inequality,
\begin{align*}
\inf_{w\in \Gamma} \abs{\E E_N(w)} &\geq  \inf_{w\in \Gamma} \abs{\frac{N}{\abs{w}}- \frac{1}{\abs{\xi-w}}}=\abs{\frac{N}{\abs{\xi}-\frac{\abs{\xi}}{N-1}} - \frac{1}{r\abs{\xi}/N}}=\frac{N}{\abs{\xi}}\lr{\frac{N-1}{N-2}-\frac{1}{r}}. 
\end{align*}
If $X_{N,\xi,\Delta}$ happens, we also have
\begin{align*}
\sup_{w\in \Gamma} \abs{\twiddle{E}_N(w)} \leq \frac{N^{1-\gamma}}{\abs{\xi}}.
\end{align*}
Hence, on the event $X_{N,\xi,\Delta}\cap B_{N,\xi}$ for which
\[\P\lr{X_{N,\xi,\Delta}\cap B_{N,\xi}} = \P\lr{X_{N,\xi,\Delta}}+O\lr{N^{-1+2\Delta}}\]
we find 
\[\inf_{w\in \Gamma}\abs{\E{E_N}(w)}>\sup_{w\in \Gamma}|\twiddle{E}_N|.\] 
We may therefore apply Rouch\'e's Theorem to conclude that $E_N$ has exactly one zero (and hence $p_{N,\xi}$ has exactly one critical point) in the interior of $\Gamma$. This is precisely the first condition in the definition of $S_{N,\xi}.$ 

To check that the small petal condition in the definition of $S_{N,\xi}$ is also satisfied when $X_{N,\xi,\Delta}$ occurs, let $A_{N,\xi}$ denote the annulus centered at $\xi$ with inner radius $r \abs{\xi}/N$ and outer radius $4r \abs{\xi} /  N $ (recall that $r$ was fixed by \eqref{E:r-def}). For simplicity, we will rotate our coordinates so that $\xi$ lies on the positive real axis and consider the three regions in $A_{N,\xi}$ (see Figure \ref{F:annulus} below):
\begin{align*}
\Omega_1 &:= \{w \in A_{N,\xi} : \Re w < \xi - 3r \xi/N \}, \\
\Omega_2 &:= \{w \in A_{N,\xi} : \xi - 2r \xi/N < \Re w < \xi + 2r\xi/N \}, \\
\Omega_3 &:= \{w \in A_{N,\xi} : \Re w > \xi + 3r \xi/N \}.
\end{align*}

\begin{figure}[h]
 \includegraphics[scale=.3]{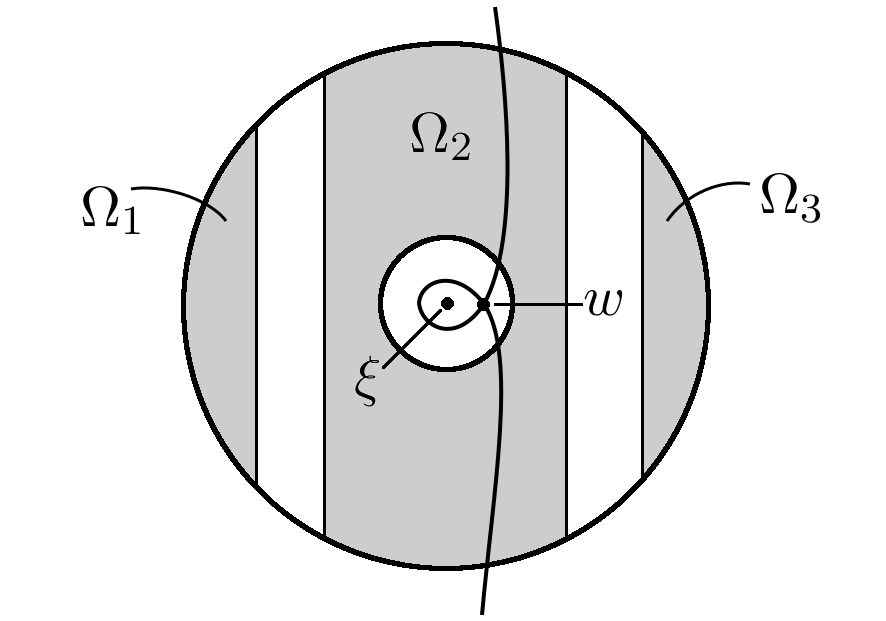}
 \caption{The lemniscate $\{|p_N(z)| = \abs{p_N(w)}\}$ passing through $w$ 
 has with high probability empty intersection with each of $\Omega_1$ and $\Omega_2$.
 This forces one of the
 corresponding petals to be contained in $D_{N,\xi}$ 
 which implies that petal is empty (it contains no critical points in its interior).}
  \label{F:annulus}
\end{figure}

We now argue that the event $X_{N,\xi,\Delta}$ in (\ref{eq:disk}) implies that the argument of $E_N(w)$ is essentially deterministic (given by the argument of $\E{E_N(w)}$ to leading order in $N$) uniformly as $w$ ranges over $A_{N,\xi}$. We parameterize points $w \in A_{N,\xi}$ by writing $w = \xi + \rho e^{i \theta}$ with $r \abs{\xi}/ N < \rho < 4r  \abs{\xi}/N$. Then (\ref{eq:expectedfield}) yields
\begin{equation}\label{eq:parentheses}
\E{E_N(w)} = \frac{N}{\abs{\xi}}\left(-\frac{1}{1 + \rho e^{i \theta}/\abs{\xi}} + \frac{1}{N \rho e^{i \theta}/\abs{\xi}} \right).
\end{equation}
The first term inside the parentheses in (\ref{eq:parentheses}) is 
$-1 + O(1/N)$ since
\[\abs{-\frac{1}{1 + \rho e^{i \theta}/\abs{\xi}} +1} \leq \frac{4r}{N-4r}=O\lr{\frac{1}{N}},\]
while the second term satisfies
\[ \frac{1}{4r} \leq \abs{  \frac{1}{N\rho e^{i \theta}/\abs{\xi}} }\leq \frac{1}{r},\]
which we summarize by writing
\[ \frac{1}{N\rho e^{i \theta}/\abs{\xi}} = F(\rho) \frac{1}{r} e^{-i\theta},\]
and taking note that 
$$\frac{1}{4} \leq F(\rho) \leq 1.$$

Therefore, Proposition \ref{prop:main} (along with the estimates above) 
shows that for each $\Delta\in (0,1/4)$ there exists $\gamma>0$ so that with probability at least $1-C_{\Delta}N^{-2\Delta}$, we have
\[E_N(\xi + \rho e^{i\theta}) = \E E_N(\xi + \rho e^{i\theta}) + \tilde{E}_N(\xi + \rho e^{i\theta}) = \frac{N}{\abs{\xi}}\lr{-1 + F(\rho)\frac{1}{r}e^{-i\theta} +O\lr{N^{-\gamma}}},\]
where the implied constant is independent of $\rho, \theta,N,\xi.$ Thus, using the definition \eqref{E:r-def} of $r$ we conclude when $X_{N,\xi,\Delta}$ occurs, we also have
\begin{equation}
\arg\lr{E_N(w)}\in \lr{\pi - \frac{1}{5}, \pi + \frac{1}{5}},\qquad \forall w\in A_{N,\xi}\label{E:avg-arg}
\end{equation}
for all $N$ sufficiently large. Let us write
\[\partial_-\Omega_2= \set{w\in A_{N,\xi}~|~\Re(w)=\xi - 2r\xi/N},\quad \partial_+\Omega_2= \set{w\in A_{N,\xi}~|~\Re(w)=\xi + 2r\xi/N}\]
for the left and right boundaries of $\Omega_2$ and set
\[\zeta = \pi/2 - \arctan 8>0.\]
Note that the angle of any line segment joining a point on $\partial_- \Omega_2$ to any point in $\Omega_1$ lies in the interval $(\pi/2+\zeta,3 \pi/2 - \zeta)$, so that it forms an acute angle with $\E{E_N (w)}$ when $X_{N,\xi,\Delta}$ happens by \eqref{E:avg-arg}. Since $E_N(w)$ has the same argument as the gradient of $|p_{N,\xi}(w)|$, this implies that the event $X_{N,\xi,\Delta}$ entails that the directional derivative of $|p_{N,\xi}(w)|$ along such a line segment is positive, and hence the value of $|p_{N,\xi}(w)|$ in $\Omega_1$ is strictly larger than its value on $\partial_-\Omega_2$. Similarly, the value of $|p_{N,\xi}(w)|$ on the right boundary $\partial_+\Omega_2$ is strictly larger than its value throughout $\Omega_3$. This implies that a level curve of $|p_{N,\xi}|$ that intersects $\Omega_2$ cannot intersect $\Omega_1$ or $\Omega_3$ unless it leaves $D_{N,\xi}$ 
through the set
$$S:=\{|w-\xi|=4r \xi/N\} \cap \Omega_2,$$
which consists of two circular arcs symmetric with respect to the real axis. 
As above, the event $X_{N,\xi, \Delta}$ ensures that the argument of $E_N$ is close to $\pi$ and hence the restriction of $|p_{N,\xi}|$ to each component of $S$ is strictly monotone. Thus, any level curve of $|p_{N,\xi}|$ can only cross each component of $S$ once on the event $X_{N,\xi,\Delta}$. 
The singular component $\Gamma_{w}$ (consisting of two petals joined at $w$)  of the lemniscate passing through the critical point $w$ that is paired to $\xi$ therefore crosses the boundary of $D_{N,\xi}$ at most twice (one crossing for each component of $S$). 
This implies whenever $X_{N,\xi,\Delta}$ occurs, one of the petals is completely contained in $D_{N,\xi}$, and therefore it must be a small petal since $D_{N,\xi}$ contains only one zero of $p_N$, see Figure \ref{F:annulus}.
This shows that $X_{N,\xi,\Delta}\cap B_{N,\xi}$ implies $S_{N,\xi}$ and yields \eqref{E:Goal1}, completing the proof of Theorem \ref{thm:poly}. \hfill $\square$

\subsection{Proof of Proposition \ref{prop:main}}\label{S:main-prop-pf}
Fix $\Delta \in (0,1/4)$ and a sequence $\xi=\xi(N)$ with $\abs{\xi}\leq N^{\Delta}$ (we remind the reader that $\abs{\xi}$ is measured 
in coordinates centered at $\infty$,
and hence in terms of the original coordinates
our assumption removes a disk of radius $N^{-\Delta}$ centered at $0$). 
In this section we explain how to modify the proof of Theorem 1 
(specifically equation (4.2)) in \cite{Hanin} to prove Proposition \ref{prop:main}. 
The argument from \cite{Hanin} was presented in several steps;
below we explain the modifications needed at each step.
\\
\noindent {\bf Step 1.} With $w_\xi = w_{\xi,N}$ 
defined as in \eqref{E:w-def}, 
we study $\twiddle{E}_N(w)$ by separately considering
\[\twiddle{E}_N(w_\xi),\qquad \text{and}\qquad \twiddle{E}_N(w)-\twiddle{E}_N(w_\xi).\]
\noindent {\bf Step 2.} To understand 
\[\twiddle{E}_N(w)-\twiddle{E}_N(w_\xi)=\sum_{j=1}^{N-1}\frac{w_\xi-w}{(w-\xi_j)(w_\xi-\xi_j)}\] 
we first fix $\delta\in (2\Delta, 1)$ 
and estimate the contribution from zeros far away from $\xi$: 
\begin{equation}
\abs{\sum_{\abs{\xi-\xi_j}>N^{-1/2+\delta/2}}\frac{w_\xi-w}{(w-\xi_j)(w_\xi-\xi_j)}}\leq K_1 N^{1-\delta+\Delta},\label{E:far-zeros}
\end{equation}
for some $K_1 > 0$, where we've used that $\abs{w-w_\xi}\leq 4r N^{-1+\Delta}$
for $w\in D_{N,\xi}$ and that $\abs{w-\xi_j},\, \abs{w_\xi-\xi_j}=\Theta(\abs{\xi-\xi_j}).$ Hence, as long as
\[\delta > 2\Delta,\]
we find that the expression in \eqref{E:far-zeros} is (deterministically) bounded above by 
\[N^{1-\ep-\Delta}\leq N^{1-\ep}/\abs{\xi}\]
for some $\ep>0$, as in the definition of the event $X_{N,\xi}$ whose probability we seek to estimate. \\

\noindent {\bf Step 3.} 
Next, we control the contribution to $\twiddle{E}_N(w)-\twiddle{E}_N(w_\xi)$ from zeros near $\xi$ by repeatedly adding and subtracting 
$\sum_{\abs{\xi-\xi_j}\leq N^{-1/2+\delta/2}}\frac{1}{(w_\xi-\xi_j)^k}$:
\begin{align}\label{E:I-II}
\abs{\sum_{\abs{\xi-\xi_j} \leq N^{-1/2+\delta/2}}
  \frac{w_\xi-w}{(w_\xi-\xi_j)(w-\xi_j)}}&\leq \abs{\sum_{\abs{\xi_j-\xi}\leq N^{-1/2+\delta/2}}
\frac{\lr{w_\xi-w}^L}{(w_\xi-\xi_j)^L (w-\xi_j)}}\\
&+\sum_{k=1}^{L-1} c^k N^{(-1+\Delta)k}\abs{\sum_{\abs{\xi-\xi_j}\leq N^{-1/2+\delta/2}}\frac{1}{(w_\xi-\xi_j)^{k+1}}}, 
\end{align}
where $c$ is an absolute constant and $L$ is any positive integer.  \\

\noindent {\bf Step 4.} We control the two terms in \eqref{E:I-II} separately. To control the term containing the sum on $k,$ we use \cite[Lem. 2]{Hanin}, which says that for every $\eta\in (0,1/2)$ there exists $K=K(\eta)>0$ so that
\[\P\lr{\sum_{j=1}^{N-1}\frac{1}{\abs{w_\xi-\xi_j}^2}>N^{2-2\eta}}\leq K \cdot N^{-1+2\eta}\log N.\]
Hence, taking $\eta=\Delta+\ep/2,$ we find
\begin{align*}
  \sum_{k=1}^{L-1} c^k N^{(-1+\Delta)k}\abs{\sum_{\abs{\xi-\xi_j}\leq N^{-1/2+\delta/2}}\frac{1}{(w_\xi-w)^{k+1}}}&\leq \sum_{k=1}^{L-1}c^k N^{(-1+\Delta)k} N^{(2-2\eta)\frac{k+1}{2}}\\
&\leq C_\ep N^{1-\Delta-\ep}\leq C_\ep N^{1-\ep}/\abs{\xi}
\end{align*}
with probability at least $N^{-1+2\Delta + \ep}.$\\

\noindent{\bf Step 5.} To bound the other term in \eqref{E:I-II}, 
we use \cite[Lem. 1]{Hanin}, which says that for each $\delta\in (0,1/2)$ with probability at least $1-C_\delta N^{-\delta}$ there are no zeros with $\abs{\xi_j}> N^{1/2+\delta/2}$ and at most $N^{2\delta}$ zeros with $\abs{\xi}>N^{1/2-\delta/2}$. This allows use to write
\[ \abs{\sum_{\abs{\xi_j}-\xi\leq N^{-1/2+\delta/2}}
\frac{\lr{w_\xi-w}^L}{(w_\xi-\xi_j)^L (w-\xi_j)}}\leq N^{1-\lr{1-2\delta+L(1-\Delta) - (L+1)(1/2 + \delta/2)}}.\]
Hence, taking $L$ sufficiently large, we find that if $2\Delta < \frac{1}{2},$ then the left hand side in the previous line can be bounded above by $N^{-1+\ep}/\abs{\xi}$ for all $\ep>0$ sufficiently small with probability at least $1-C_{\Delta} N^{-1+\Delta}.$ 
This completes the proof of Proposition \ref{prop:main}.

\section{Random perturbation of a Chebyshev polynomial}\label{S:Cheby}

In light of the results of the previous section, 
one may wonder whether there are any natural models
of random polynomials that typically have some 
positive portion of the nodes 
in the corresponding tree having two children
(thus resembling the previously established combinatorial baseline).
As a possible candidate for such a model the authors considered 
random linear combinations of Chebyshev polynomials.
More specifically, the class of polynomials considered were of the form $p(z) =  \sum_{k=0}^n a_k T_k(z)$, where $T_k$ is the Chebyshev polynomial (of the first kind) of degree $k$, and the coefficients are chosen independently with $a_n \sim N(0,1)$. 
Linear combinations of orthogonal polynomials have been studied previously,
including several varieties of Jacobi orthogonal polynomials \cite{LPX}.
The important property of Chebyshev polynomials (leading us to choose 
those as a basis) is that they each have critical values all with the same modulus.
Figure \ref{fig:Cheb} shows the family of singular level sets for such a polynomial. 
The lemniscates for this type of polynomial appear to exhibit a 
rich nesting structure.
However, these polynomials are typically not lemniscate generic
(due to complex conjugate pairs of critical points sharing the same critical value).
Consequently, 
this model seems worthy of further investigation,
but this will require first understanding an appropriate class of lemniscate trees.\\

We will investigate a modified (less organic, but more tractable)
version of this model where the top degree Chebyshev polynomial
gets most of the weight.
Specifically, we consider randomly perturbed Chebyshev polynomials of the form 
$T_n(z) + \frac{1}{n} \sum_{k=0}^{n-1} b_k T_k(z)$, 
where the coefficients $b_k$ are randomly and independently chosen to be 
$1$ or $-1$ with equal probability. 
These polynomials have all real roots and real critical points, 
which enables us to easily determine the corresponding lemniscate 
tree by the process described below.\\

Suppose that $p$ is a lemniscate generic polynomial with real zeros and critical points. 
We construct a permutation as follows: 
we label the critical points with the integers 1 through $\text{deg}(p)-1$, starting with 1 for the critical point with largest critical value in magnitude, 2 for the critical point with second largest critical value in magnitude, and so on. Reading the labels from left to right gives a permutation of the numbers 1 through $\text{deg}(p)-1$. Now, it is well known that the permutations on $n$ letters are in one-to-one correspondence with the \textit{increasing binary trees} of size $n$ (see \cite[p. 143]{Flajolet} for example). These are plane, labeled, rooted trees in which every vertex has at most two children, where each child has a left or right orientation (even when it is the unique child of its parent), such that the labels along any path directed away from the root are increasing. We then construct the increasing binary tree corresponding to the permutation obtained from the polynomial. By "forgetting" its embedding in the plane we obtain the lemniscate tree associated to the singular level sets of the polynomial. One can even determine the number of nodes of outdegree 2 directly from the permutation by counting the number of descents which are immediately followed by ascents.\\

\begin{figure}
\centering
\begin{subfigure}{.5\textwidth}
  \centering
  \includegraphics[width=0.8\linewidth]{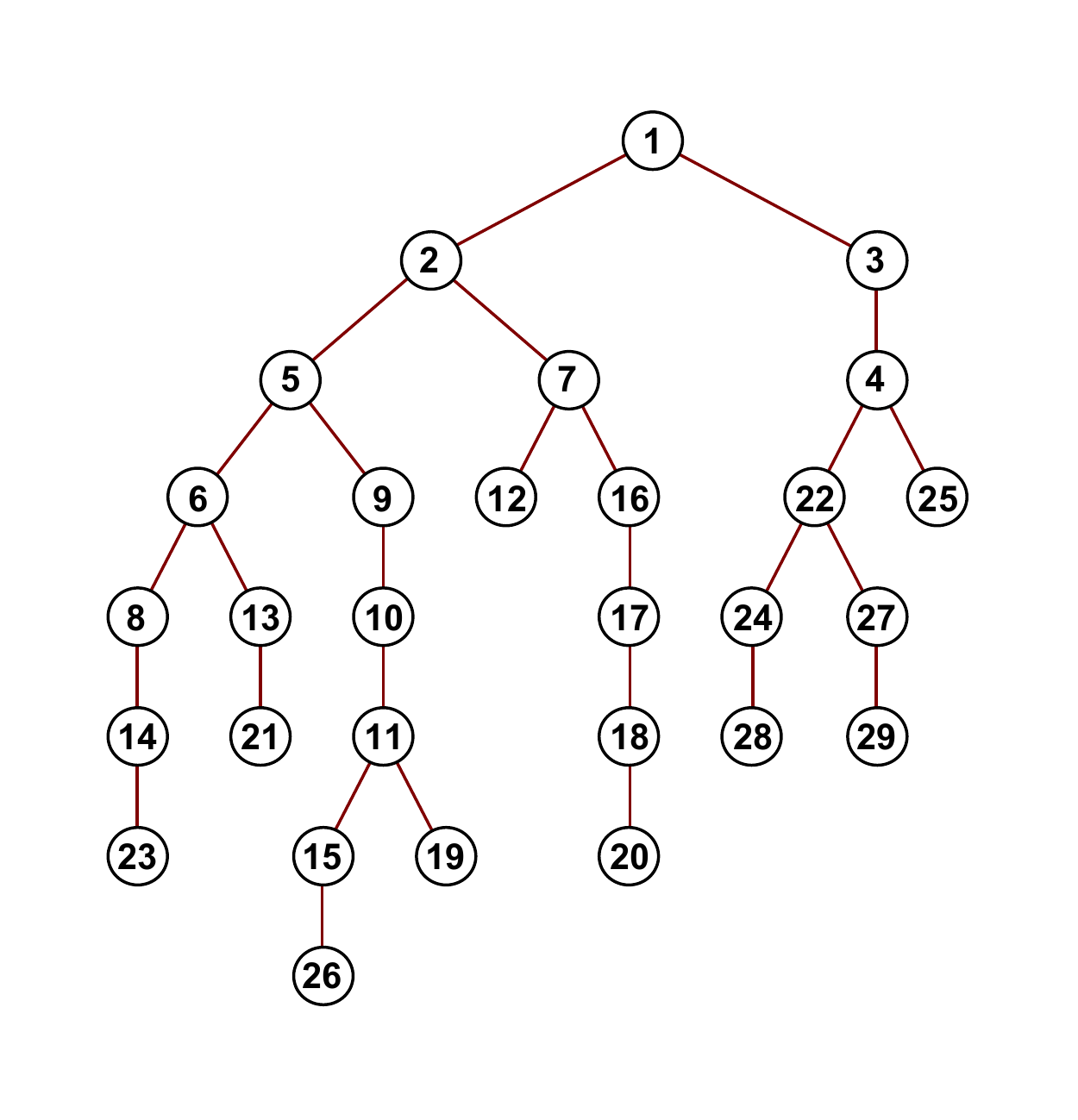}
\end{subfigure}%
\begin{subfigure}{.5\textwidth}
  \centering
  \includegraphics[width=0.8\linewidth]{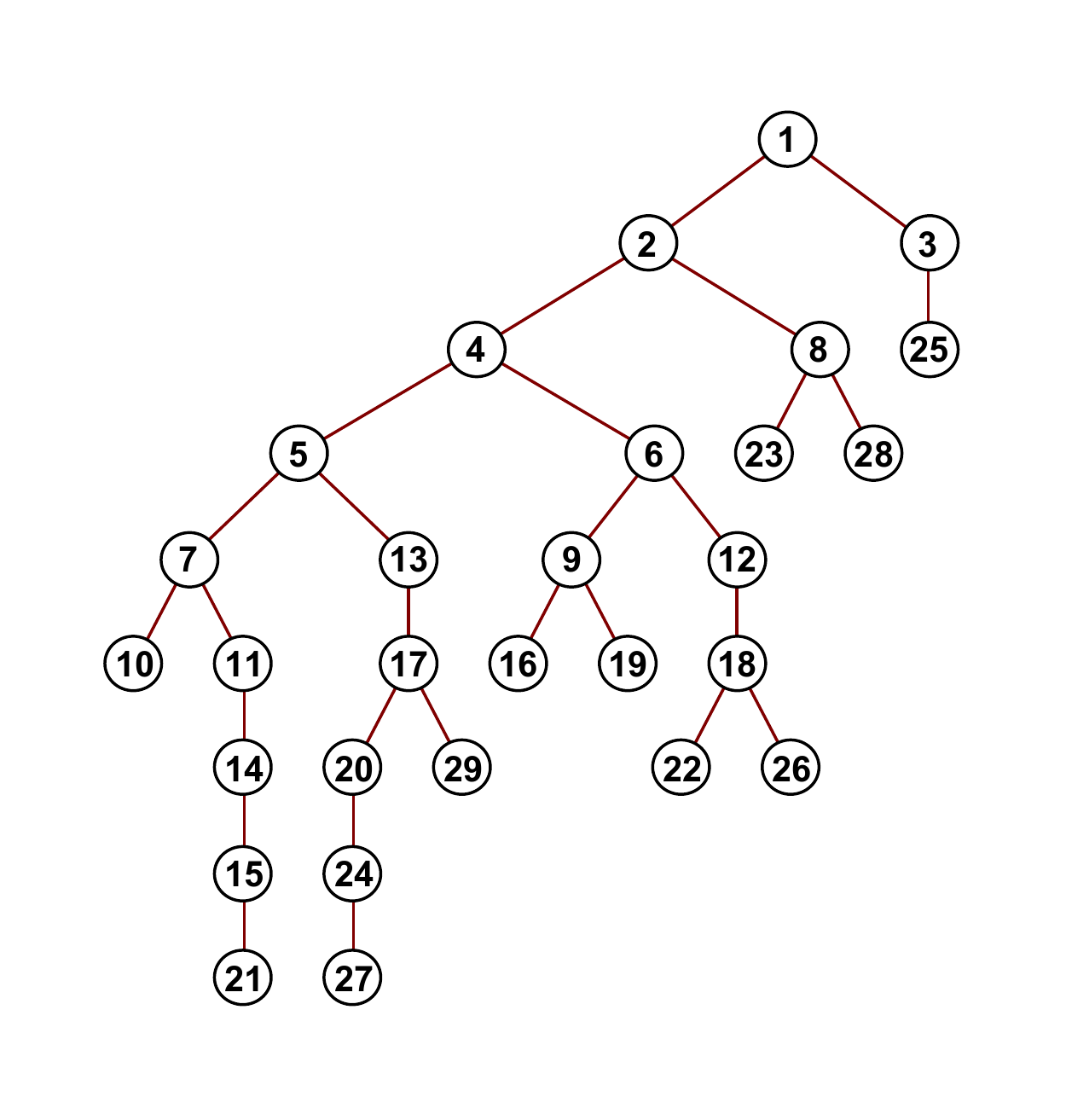}
\end{subfigure}
\caption{Lemniscate trees for randomly perturbed Chebyshev polynomials of degree 30.}
\label{fig:trees}
\end{figure}

We apply this procedure to a number of polynomials in the following computer experiment: Table \ref{tab:meanN2} gives the average value of $N_2$ computed for a sample of $100$ randomly perturbed Chebyshev polynomials of the same degree $n$ for a number of different values of $n$ ranging from $10$ to $200$. Linear regression yields a best fit line with equation $N_2 = 0.3338n -0.90803$ with $R^2 = 0.9999$, indicating that one should expect 
approximately a third of the vertices in the lemniscate tree for a 
perturbed Chebyshev polynomial to have outdegree two.
This agrees with a heuristic of
ignoring correlations in the randomly perturbed heights of critical values
in the perturbed Chebyshev polynomial,
which corresponds to the induced random permutation being sampled 
uniformly from the combinatorial class of permutations
(it is known \cite{Bergeron} that the average number of nodes of outdegree two
in a random permutation tree is asymptotically a third of the nodes).
Figure \ref{fig:trees} shows the lemniscate trees corresponding to two 
randomly perturbed Chebyshev polynomials of degree 30.

\begin{table}[hbt]
\caption{Average value of $N_2$ vs. degree}
\label{tab:meanN2}
\begin{center}
\begin{tabular}{c|cccccccccc}

$n$ & 10 & 20 & 30 & 40 & 50 & 60 & 70 & 80 & 90 & 100\\ \hline \noalign{\smallskip}

mean $N_2$ & 2.55 & 5.79 & 9.23 & 12.53 & 15.47 & 19.01 & 22.27 & 25.63 & 29.1 & 32.64\\   \noalign{\medskip}

$n$ & 110 & 120 & 130 & 140 & 150 & 160 & 170 & 180 & 190 & 200\\ \hline \noalign{\smallskip}

mean $N_2$ & 35.77 & 39.39 & 42.42 & 46.09 & 49.06 & 52.73 & 55.86 & 59.08 & 62.44 & 65.72

\end{tabular}
\end{center}
\end{table}

%\begin{center}
%\begin{figure}
%\caption{}
%\label{line}
%\includegraphics{N2vsDegree.png}
%\end{figure}
%\end{center}

\begin{remark}
Random matrix theory gives rise to a
more natural model of random polynomials 
that may yet exhibit a similar outcome
as the perturbed Chebyshev model.
Namely, consider the characteristic polynomial
$p(z) = \det(M-z I)$ of a random matrix 
$M$ sampled from the so-called Jacobi ensemble \cite{DE}
(with parameters chosen in order that the associated Jacobi orthogonal polynomials
are Chebyshev polynomials).
We expect that $p(z)$ has a lemniscate tree with, on average,
approximately one third of its nodes of outdegree two.
\end{remark}

%\comm{The idea behind using the Chebyshev model is that the critical values are at comparable heights, but is there an electrostatic interpretation?  The zeros of a single Chebyshev polynomial are at Fekete points and solve a minimum energy problem in potential theory, but I do not see why this causes them to shield themselves from the effect of a high-order pole at infinity.}

\section{Acknowledgements}
The authors would like to thank Alexandra Milbrand for 
creating the picture of the landscape of level sets in figure \ref{fig:illustration}.
We would also like to thank Guilherme Silva for pointing out the reference
\cite{DE} to a relevant ensemble of random matrices.

\vspace{0.1in}
 
{\em 

Email: mepstein2012@fau.edu

Email: bhanin@math.tamu.edu

Email: elundber@fau.edu}

\vspace{0.1in}

\end{document}